\documentclass[12pt,notitlepage]{amsart} 
\usepackage{amssymb,epic,latexsym} 
\newlength{\fw}
\newlength{\jmr}
\newlength{\jfc}
\newlength{\bernd}
\newlength{\ioan}
\newlength{\agk}
\newlength{\gkz}
\newtheorem{lemma}{Lemma}

\newtheorem{dfn}{Definition}
\newtheorem{cordfn}[dfn]{Corollary and Definition}
\newtheorem{van}{Vanishing Theorem for Resultants} 

\newtheorem{main}{Main Theorem} 

\newtheorem{thm}{Theorem}
\newtheorem{rem}{Remark}

\newcommand{\relint}{\mathrm{RelInt}}

\newcommand{\cA}{\mathcal{A}}
\newcommand{\cD}{\mathcal{D}}
\newcommand{\cE}{\mathcal{E}}

\newcommand{\cH}{\mathcal{H}}

\newcommand{\cP}{\mathcal{P}}

\newcommand{\supp}{\mathrm{Supp}}
\newcommand{\conv}{\mathrm{Conv}}

\newcommand{\thth}{{\underline{\mathrm{th}}}}
\newcommand{\rd}{{\underline{\mathrm{rd}}}}

\newcommand{\Pro}{\mathbb{P}}
\newcommand{\Pa}{\Pro^{|A|-1}_K}

\newcommand{\R}{\mathbb{R}}
\newcommand{\C}{\mathbb{C}}
\newcommand{\N}{\mathbb{N}}
\newcommand{\Z}{\mathbb{Z}}
 
\newcommand{\fii}{\varphi}

\newcommand{\ch}{\mathrm{char}}

\newcommand{\divisor}{\mathrm{Div}}

\newcommand{\choo}{\mathrm{Ch}}
\newcommand{\fan}{\mathrm{Fan}}

\newcommand{\res}{\mathrm{Res}}

\newcommand{\Sn}{\mathcal{S}^{n-1}}
\newcommand{\Zn}{\Z^n}

\newcommand{\Rn}{\R^n}

\newcommand{\Ks}{K^*}
\newcommand{\Kn}{K^n}

\newcommand{\Rns}{\Rn\!\setminus\!\{\bO\}} 

\newcommand{\Ksn}{{(K^*)}^n}

\newcommand{\cM}{\mathcal{M}}

\newcommand{\wP}{{\widetilde{P}}} 
\newcommand{\bZ}{{\bar{Z}}} 
\newcommand{\hZ}{\bar{Z}^\vee_+} 
\newcommand{\dZ}{\bar{Z}^\vee}

\newcommand{\cF}{\mathcal{F}}

\newcommand{\cT}{\mathcal{T}}
\newcommand{\bO}{\mathbf{O}}

\begin{document}
\title{Toric Generalized Characteristic Polynomials} 

\author{J. Maurice Rojas}
\thanks{This research was completed at MIT and was 
partially funded by an N.S.F. Mathematical Sciences Postdoctoral Fellowship.
Research at MSRI is supported in part by NSF grant DMS-9022140.}

\address{\hskip-\parindent
Massachusetts Institute of Technology\\
	 Mathematics Department\\
         77 Mass. Ave.\\
         Cambridge, MA \ 02139, U.S.A. }

\date{\today}

\begin{abstract}
We illustrate an efficient new method for handling polynomial systems with 
degenerate solution sets. In particular, a corollary of our techniques 
is a new algorithm to find an isolated point in every excess component 
of the zero set (over an algebraically closed field) of any $n$ by $n$ 
system of polynomial equations. Since we use the sparse resultant, we 
thus obtain complexity bounds (for converting any input 
polynomial system into a multilinear factorization problem) which are close to 
cubic in the degree of the underlying variety --- significantly better 
than previous bounds which were pseudo-polynomial in the classical B\'ezout 
bound. By carefully taking into account the underlying toric 
geometry, we are also able to improve the reliability 
of certain sparse resultant based algorithms for polynomial system 
solving.
\end{abstract} 

\maketitle

\section{Introduction}
\label{sec:intro}
The rebirth of resultants, especially through the {\bf 
toric}\footnote{Other commonly used prefixes for this more modern 
generalization include: sparse, mixed, $\cA$-, and $(\cA_1,\ldots,\cA_k)$-.} 
resultant 
\cite{gkz94}, has begun to provide a much needed alternative to 
Gr\"obner basis methods for solving polynomial systems. Continuing 
this philosophy, we will present a new, fast and reliable, resultant-based 
method for handling certain degenerate polynomial systems. Simply put, 
we refine and generalize the {\bf generalized characteristic polynomial} 
(GCP) \cite{gcp,mikebezout} to take {\bf sparsity}\/ into account. Furthermore, 
we introduce the concept of a {\bf twisted Chow form} in order to 
completely avoid any degeneracies within our algorithm. 

The importance of dealing with degenerate polynomial systems has 
been observed in earlier work on quantifier elimination over 
algebraically closed fields \cite{chigo,renegar,pspace}: 
Many reasonable algorithms fail catastrophically when 
presented with an $n$ by $n$ system having positive-dimensional zero 
set. When such is the case, it is of considerable benefit to the 
user to at least be given some sort of description of the zero-dimensional  
part of the zero set. This was a benefit of Canny's original 
GCP, but he remarked \cite[pg.\ 242]{gcp} ``For large, 
dense problems however, the resultant and GCP methods should be faster 
[than Gr\"obner algorithms].'' Our construction, 
the {\bf toric GCP}, promises to be much more competitive when applied 
to {\bf sparse} systems in such a comparison.

\begin{rem}
It should be emphasized that perturbation methods for degenerate 
systems (such as the toric GCP) are of the greatest importance 
when working with {\bf exact} arithmetic. However, floating point 
polynomial system solving also benefits from a complete and rigourous 
understanding of the potential degeneracies within exact arithmetic, e.g., 
\cite{smale87}. In any case, it is frequently the case in many applications 
that ``real life'' happens to land in a measure $0$ exception which breaks an 
algorithm.
\end{rem}

In what follows, we will frequently use multi-index notation in 
order to precisely control which monomial terms are allowed to appear in 
our polynomial systems.\footnote{This is what we mean by sparsity, 
as opposed to working solely in terms of polynomial degrees. Note 
that our notion overlaps the case of having very few monomial terms.} 
This notation, and in particular, 
{\bf supports} and {\bf Newton polytopes}, are amply detailed 
in earlier works of Emiris \cite{emiphd}, Huber \cite{birkphd}, Gelfand 
et.\ al.\ \cite{gkz90,gkz94}, Rojas \cite{convexapp,toricint}, 
Sturmfels \cite{sparseelim}, and Verschelde \cite{janphd}, to mention 
but a few authors and references. 

Since what we will present is at heart a perturbation method, we 
will first need the following definition to construct certain 
polynomial systems in ``general position.'' 
\begin{dfn}
\cite{convexapp,rojaswang}
Given $n$-tuples $D\!:=\!(D_1,\ldots,D_n)$ and $E\!:=\!(E_1,\ldots,E_n)$ 
of nonempty finite subsets of $\Zn$, we say that $D$ 
{\bf fills} $E$ iff (0) $D_i\!\subseteq\!E_i$ for all 
$i\!\in\![n]$ and (1) $\cM(D)\!=\!\cM(E)$. 
An {\bf irreducible} fill is then simply a fill which is minimal with respect 
to $n$-tuple containment. If $D$ or $E$ is instead an $n$-tuple of 
polytopes in $\Rn$, then we will use the same definition. 
\end{dfn}
In the above, $\cM(E)$ denotes the {\bf mixed volume} 
\cite{buza,schneider,isawres,mvcomplex}
of the convex hulls of the $E_i$, and we use $[j]$ for the set of 
integers $\{1,\ldots,j\}$.

Our construction is summarized in the following definition and 
main result. (Henceforth, all of our polynomials and roots are to be 
considered over an algebraically closed field $K$.) 
\begin{dfn} 
\label{dfn:main1}
Suppose $F$ is an $n\!\times\!n$ polynomial system with support contained 
in $E$ and $g(x)\!:=\!\sum_{e\in A} u_e x^e$, where the $u_e$ are algebraically 
independent indeterminates and $A\!\subset\!\Zn$ is nonempty and finite. 
Assume further that $\cM(E)\!>\!0$ and $A$ has at least two elements. 
Letting $D$ be an irreducible fill of $E$, $F^*\!:=\!(\sum_{e\in D_i} x^e \; | 
\; i\!\in\![n])$, and $u\!:=\!(u_e \; | \; e\!\in\!A)$, define 
${\boldsymbol \choo_A(u)}\!:=\!\res_{(E,A)}(F,g)$ and 
${\boldsymbol \cH(u;s)}\!:=\!\res_{(E,A)}(F-sF^*,g)$, where $s$ is a new 
indeterminate. We call $\cH$ a {\bf toric} {\bf generalized 
characteristic polynomial for $\mathbf{(F,A)}$.} 
\end{dfn}
\noindent
In the above, $\res_\star(\cdot)$ denotes the toric resultant. Recall 
also that to any $n$-dimensional rational polytope $P\!\subset\!\Rn$ one can 
associate its corresponding {\bf toric variety} $\cT_P$ 
\cite{dannie,ksz92,tfulton,toricint}.
\begin{main}
Following the notation of definition \ref{dfn:main1}, 
consider $\cH$ as a polynomial in $s$ with coefficients in $K[u]$ and let 
${\boldsymbol \cF_A(u)}$ be the coefficient of the lowest degree term of 
$\cH$. Also let $P\!:=\!\sum^n_{i=1}\conv(E_i)$, $\wP\!:=\!P+\conv(A)$, and 
$\fii : \cT_\wP \longrightarrow \Pa$ the natural morphism defined 
by  $x \mapsto [x^e \; | \; e\!\in\!A]$ (cf.\ section \ref{sec:chowtor}). Then 
$\cF_A$ is a homogeneous polynomial, of degree $\cM(E)$, with the following 
properties:
\begin{enumerate}
\addtocounter{enumi}{-1}
\item{The constant term of $\cH(s)$ is precisely $\choo_A$. In particular, 
letting $Z$ be the zero set of $F$ in $\cT_\wP$, $\choo_A\!\equiv\!0 
\Longleftrightarrow \fii(Z)$ is positive-dimensional. } 
\item{If $\zeta\!\in\!\cT_\wP$ is an isolated root of $F$ then 
$\cF_A$ is divisible by $\sum_{e\in A} c_e u_e$, where 
$[c_e \; | \; e\!\in\!A]\!=\!\fii(\zeta)$.}  
\item{If a nonzero linear form $\sum_{e\in A} c_e u_e$ divides 
$\cF_A$ then $[c_e \; | \; e\!\in\!A]\!=\!\fii(\zeta)$ for {\bf some} root 
$\zeta\!\in\!\cT_\wP$ of $F$.}
\item{$\cF_A$ splits completely into linear factors. In particular, under the 
correspondence of (2), we can explicitly find at least one root of $F$ within 
every positive-dimensional irreducible component of the zero set of $F$ in 
$\cT_\wP$.}
\end{enumerate} 
We call $\cF_A$ a {\bf toric perturbation} of $\choo_A$. 
\end{main} 
A complete criterion for finding fills is given in theorem \ref{thm:me} 
of section \ref{sec:fill} and some simple examples of filling 
(and our Main Theorem) appear in section \ref{sec:ex}. Compatibility is 
defined in section \ref{sec:chowtor}. There we 
also detail what is meant by the zero set of $F$ in a toric variety.  
\begin{rem}
Alternatively, if one wants to avoid filling, one can substitute for $F^*$ 
any polynomial system 
(with support contained in $E$) {\bf known} to have exactly $\cM(E)$ 
roots (counting multiplicities) in $\Ksn$. So one can also view filling as a 
deterministic way of explicitly constructing a ``generic'' polynomial system. 
\end{rem}

So how does one actually use the above theorem? One simple example 
is the {\bf sparse $\mathbf{u}$-resultant} \cite{emiphd}, which is simply a 
variant of the classical $u$-resultant \cite{vdv}. It can be defined simply
as $\choo_A(u)$ where $A$ is the vertex set of the standard $n$-simplex in 
$\Rn$. One useful (and easily verified) property of the sparse $u$-resultant 
is that $F$ has a root $\zeta\!:=\!(\zeta_1,\ldots,\zeta_n)\!\in\!\Ksn 
\Longrightarrow u_0+\zeta_1u_1+\cdots+\zeta_nu_n$ divides $\choo_A(u)$. So, 
assuming $\choo_A(u)$ is not identically $0$ (and that one has good software 
for toric resultants and multivariate factoring), one can find the 
isolated roots of $F$ simply by factoring $\choo_A(u)$. Thus the degenerate 
instances $\choo_A(u)\!\equiv\!0$ obstruct this reduction to factoring 
and our Main Theorem allows us to avoid this problem: Simply use 
$\cF_A(u)$ instead of $\choo_A(u)$. 
\begin{rem}
An interesting ``failure'' for the classical $u$-resultant 
is the case where $F$ has only finitely many roots in affine space, 
while having infinitely many in projective space. This is where 
the classical GCP is especially handy. For the toric GCP, the 
toric variety $\cT_P$ plays the role of projective space. 
This is part of the philosophy of sparse elimination theory: 
working in a well chosen toric compactification (depending 
on $F$) leads to better algorithms than if one were to work only in 
projective space. In particular, {\bf lifting} to $\cT_\wP$ helps 
us detect precisely when the sparse $u$-resultant is identicaly $0$. 
\end{rem}

However, there is a stickier subtlety which occurs with 
the sparse $u$-resultant: It is possible for $F$ to 
have only finitely many roots within $\cT_P$ with 
$\choo_A(u)$ {\bf still} vanishing identically (cf.\ section 
\ref{sub:stranger}). This is the motivation for {\bf twisted Chow forms}, 
which are defined below, and pursued at much greater length in \cite{twist}.
\begin{cordfn}[cf.\ section \ref{sec:chowtor}]
\label{dfn:twist} 
Following the notation of our Main Theorem, pick $A$ 
to be the vertices of a product of simplices with which $P$ is compatible. 
We then call $\choo_A(u)$ the {\bf twisted Chow form} of the zero set of $F$ 
{\bf with respect to $A$}. Furthermore, a twisted Chow form does {\bf not} 
vanish identically if $F$ has only finitely many roots in $\cT_P$. \qed 
\end{cordfn}
In particular, since our last construction implies that $\cT_{\conv(A)}$ 
is a product of {\bf twisted} projective spaces \cite{tfulton}, 
the coefficients of a twisted Chow form are actually multisymmetric 
functions \cite{introchow,rio} of projections (by $\fii$) of roots of $F$ 
(in $\cT_P$) onto this product. So twisted Chow forms generalize the 
$u$-resultant, the sparse $u$-resultant, and (suitably extended \cite{twist}) 
the Chow form of a projective variety \cite{introchow}. Moreover, when 
combined with the {\bf Smith normal form} \cite{unimod1,unimod2,polyhomo} and 
the toric GCP, twisted Chow forms give a considerably more reliable method for 
solving polynomial systems than the sparse $u$-resultant \cite{twist}. 

An important difference to note is that our present toric GCP is 
primarily suited for finding roots in $\Ksn$, while the original GCP is mainly 
suited (in a ``non-sparse'' way) for affine space. To {\bf completely} 
generalize and improve the GCP in affine space, it is necessary 
to use the {\bf affine}  sparse resultant and this is pursued further in 
\cite{aff}. For instance, by replacing the sparse resultant with the affine 
sparse resultant, and using {\bf $\Kn$-counting} \cite{toricint} instead of 
filling, we can actually recover Canny's GCP in the dense case. 

We close this introduction with a word on the computational complexity 
of computing the toric GCP. Neglecting  preprocessing (finding an 
irreducible fill and finding a mixed subdivision in order to 
set up the toric resultant matrices \cite{isawres}), recent work of 
Emiris, Morrain, and Pan \cite{emipan,moripan} suggests that it is possible
to find the perturbation $\cF_A$ within a number of arithmetic 
steps\footnote{Although the 
possibility of ``near-quadratic'' algorithms for polynomial system solving 
is discussed in \cite{emipan,moripan}, neither paper discusses the case where 
the zero set of $F$ is positive-dimensional, which is our main concern 
here.} which is close to {\bf cubic} in $\cM(E)$. (Indeed, Canny has 
pointed out \cite{gcp} that the original GCP can be computed in time 
close to cubic in the B\'ezout bound, which would be a special case 
of time cubic in $\cM(E)$.) Since $\cM(E)$ is much smaller than 
the B\'ezout number for most polynomial systems \cite{convexapp,toricint}, 
these preliminary results suggest that the toric 
GCP has considerable potential for practical applications. 

Let us now illustrate some of our theory.

\section{Examples}  
\label{sec:ex} 
We begin with two small examples of filling. We then see 
applications of the toric GCP to some degenerate $2 \times 2$ and $3 \times
3$ polynomial systems. Finally, we see a brief comparison of the toric 
GCP to the original GCP. 

\subsection{Filling Squares and Cubes}
For our first example, consider the pair of rectangles 
$\cP\!:=\!([0,a]\!\times\![0,b],[0,c]\!\times\![0,d])$ where
$a$, $b$, $c$, and $d$ are positive integers. Then it 
is easily verified (via theorem \ref{thm:me} of section 
\ref{sec:fill}) that the pair $$D\!:=\!(\{(0,0),(a,b)\},\{(0,d),(c,0)\})$$
fills $\cP$. In this case, the mixed area of both pairs is easily 
checked to be $ad+bc$. Note also that $D$ is a pair of oppositely 
slanting diagonals of our initial pair of rectangles (modulo 
taking convex hulls). Finally, it is easily checked that $D$ is 
indeed irreducible, since the removal of any point of $D$ results in a 
mixed area of $0$. 

For our second example, let $\cP$ instead be a triple of 
standard cubes (so that the vertex set of each cube is simply 
$\{0,1\}^3$). Then, using the criterion from theorem \ref{thm:me} 
once again, it is easily verified that the triple 
\[ D\!:=\!(\{(0,0,0),
(1,1,1)\},\{(1,0,0),(0,1,0),(0,0,1)\},\{(1,1,0),(1,0,1),(0,1,1)\}) \] 
fills $\cP$. Also, it is easily checked that the mixed volume of  
both triples is $6$. Note that the triple $D$ 
consists of a body diagonal and two oppositely oriented (but  
parallel) sub-triangles of the unit cube, modulo taking convex hulls. 
Finally, note that this $D$ is irreducible as well by theorem \ref{thm:me}. 
(This is also easily checked by one of the publically available software 
packages for mixed volume by Emiris, Huber, or Verschelde.)

\subsection{A Degenerate $2\times 2$ System} 
Consider the bivariate polynomial system 
$F\!:=\!(1+2x-2x^2y-5xy+x^2+3x^3y,2+6x-6x^2y-11xy+4x^2+5x^3y)$. 
Letting $E$ be the support of $F$, the reader can easily verify that 
$\cM(E)\!=\!4$, and that the only roots of $F$ are the points 
$\{(1,1),(\frac{1}{7},\frac{7}{4})\}$ and the line $\{-1\}\!\times\!K$ 
(assuming $\ch K\!\neq\!2$).\footnote{When $\ch K\!=\!2$, the 
second isolated root becomes an isolated root lying on the $x$-axis.} So
it would appear that the $u$-resultant (and even the sparse 
$u$-resultant) will vanish identically and not give us any useful 
information about any of these roots. Let us see if perturbing 
the sparse $u$-resultant helps...

Note that by theorem \ref{thm:me}, $D\!:=\!(\{\bO,(3,1)\},\{(1,1),(2,0)\})$ is 
an irreducible fill of $E$. So applying our Main Theorem with this $D$ (and 
$A\!:=\!\{(0,0),(1,0),(0,1)\}$) one can compute with the use of {\tt Maple} 
that $\cH(u;s)$ is precisely 
\scriptsize
\[ 
(u_2^4-u_0^4+u_1^4+6u_1^2u_2^2-4u_1u_2^3-4u_1^3u_2) \mathbf{s^8} 
\]
\[+(-20u_2u_0^3-20u_2^3u_0- 4u_1u_0^3+36u_1^2u_2^2-19u_0^4-24u_2^4 
+6u_0^2u_1u_2 \]
\[ +36u_1u_2^3+36u_1^4-12u_0u_1
^2u_2-9u_1^2u_0^2+3u_2^2u_0^2+36u_0u_1u_2^2-4u_0u_1^3-84u_1^3u_2)\mathbf{s^7}\]
\[ +(-170u_2
u_0^3-394u_1^3u_2-98u_1u_0^3-98u_0^2u_1u_2-20u_0^4+220u_2^4+370u_2^3u_0\] 
\[
+14u_0u_1
u_2^2-110u_0u_1^3-226u_1^2u_0^2-354u_1^2u_2^2+454u_1^4-274u_0u_1^2u_2+
74u_1u_2^3)\mathbf{s^6}\]
\[ +(-1612u_0^2u_1u_2+1008u_2u_0^3+903u_0^4-624u_1u_0^3-2632u_2^3u_0-2104u_0
u_1^2u_2-970u_2^4\]
\[ -1010u_1u_2^3 +418u_1^3u_2 -2104u_0u_1u_2^2-642u_1^2u_2^2-1547
u_1^2u_0^2-936u_0u_1^3-1557u_2^2u_0^2+2204u_1^4)\mathbf{s^5}\]
\[ +(538u_0^2u_1u_2+1271u_0^4+
12253u_2^2u_0^2+6972u_2u_0^3+1929u_1^4-3075u_1^2u_2^2+654u_0u_1u_2^2\] 
\[+50u_1u_0^3 +2156u_2^4-960u_1^2u_0^2-2290u_0u_1^3+132u_1u_2^3-
5344u_0u_1^2u_2-1142u_1^3u_2+8708u_2^3u_0)\mathbf{s^4} \] 
\[ +(4384u_1u_0^3-24988u_2^2u_0^2-1582u_1^3u_2-6756u_0^4+10884u_0
u_1u_2^2+3802u_1u_2^3+15438u_0u_1^2u_2\] 
\[ 
+1024u_0u_1^3+8324u_1^2u_0^2-12826u_2^3
u_0+11270u_0^2u_1u_2-6976u_1^4+7164u_1^2u_2^2-21326u_2u_0^3-2408u_2^4)
\mathbf{s^3}\]
\[ +(3436
u_1^3u_2+3800u_0u_1^3+7756u_2^3u_0-3886u_1u_2^3+1225u_2^4+17059u_2^2u_0^2-5984
u_1^2u_0^2\]
\[ +15708u_2u_0^3-12232u_0u_1u_2^2+5180u_0^4-2091u_1^2u_2^2-6828u_0u_1^2
u_2+1316u_1^4-12700u_0^2u_1u_2-4312u_1u_0^3)\mathbf{s^2} \] 
\normalsize
\[ +(-196u_2^4-448u_0^4-768u_1^3u_2 +512u_1^2u_0^2-1288u_2^3u_0\]
\[ -2436u_2^2u_0^2+1920u_0u_1u_2^2 +1536u_0^2u_1u_2+1024 u_0u_1^2u_2\] 
\[ +384u_1u_0^3-1792u_2u_0^3-384u_0u_1^3+768u_1u_2^3+260u_1^2u_2^2-64
u_1^4)\mathbf{\underline{s}.} \]  
(Sparse resultants of this size are quite amenable with the 
aid of {\tt Maple}, following the technique of a similar computation in 
\cite{rio}.) So our $\cF_A$ is just the coefficient of $s$ or $s^2$ in this 
polynomial, according as $\ch K\!\neq\!2$ or $\ch K\!=\!2$. To simplify our 
discussion, let us henceforth assume the former possibility. 

Factoring with {\tt Maple}, we obtain that $\cF_A$ can be written as
follows:
\[ -4(u_0+u_1+u_2)(28u_0+4u_1+49u_2)(u_0-u_1+u_2)(4u_0-4u_1+u_2) 
\] 
In particular, given any factor above, the ratio of the coefficients 
of $u_i$ and $u_0$ is precisely the $i^\thth$ coordinate of a corresponding 
root of $F$. Thus the first two factors correspond precisely to the 
two isolated roots we already know.  As for the last two factors, 
note that they both give isolated points lying on the aforementioned 
line  $\{-1\}\!\times\! K$. This can be interpreted as assigning 
an {\bf excess intersection multiplicity} of $2$ to the line, so that
the sum of all intersection numbers (of the irreducible components of 
the zero set of $F$ in $\cT_\wP$) is $4$.

Note that the original GCP could have been used above but 
would have resulted in a $u$-form of degree $16$ (the product 
of the degrees of $f_1$ and $f_2$). Also, the corresponding 
version of $\cH(\cdot)$ is significantly larger, having 672 terms, 
compared to 110 for our above toric GCP. 

\subsection{Stranger Degeneracies}
\label{sub:stranger}
Here we give two examples showing how the sparse $u$-resultant 
can fail to find roots in $\cT_P$, even with the benefit of the toric GCP, 
unless some other construction (such as a twisted Chow form) is used.  

First consider the parameterized bivariate system 
$F\!:=\!(a_1y+a_2x+a_3xy,b_1y+b_2x+b_3xy)$. Note that 
the mixed volume bound for this system is $1$. The sparse 
$u$-resultant for this system is also easily found  
(using the same techniques as in our last example) to 
be: 
\begin{eqnarray*}
& & (a_3^2b_2b_1+b_3^2a_2a_1-b_3a_2a_3b_1-b_3a_3b_2a_1)u_0 \\
&+& (b_1^2a_2a_3-b_1b_3a_2a_1-b_2 a_1a_3b_1+b_2b_3a_1^2)u_1 \\
&+& (b_1b_3a_2^2-b_2a_2a_3b_1+b_2^2a_3a_1-b_2b_3a_2a_1)u_2 
\end{eqnarray*}
In particular, we see that when 
\[(a_1,a_2,a_3,b_1,b_2,b_3)\!=\!(0,1,2,0,1,3).\] 
the sparse $u$-resultant vanishes identically.  For this specialization, it is 
also easy to see that $F$ has only {\bf one} root in $\cT_P$, and this root 
lies at the point of $\cT_P$ corresponding to the vertex $(0,1)$ of $P$ (cf.\ 
lemma \ref{lemma:once} from section \ref{sec:chowtor}), following 
the notation of our Main Theorem. In fact, 
for this $P$, $\cT_P\!\cong\!\Pro^2_K$ and the vertex $(0,1)$ 
corresponds precisely to the point $[x:y:z]\!=\![1:0:0]$ in this particular 
copy of $\Pro^2_K$ ($z$ denoting an extra variable for homogenizing).  

There are two ways of viewing this degeneracy ($\choo_A(u)\!\equiv\!0$ while 
$F$ has only finitely many roots in $\cT_P$) of the sparse $u$-resultant. 
The first is pragmatic: One should not count this as a deficiency of the 
sparse $u$-resultant because we've cheated and set two of 
the coefficients of $F$ to $0$, thus changing the supports. 
(Our next example avoids this trick.) The second point of view 
is more geometric: By our Main Theorem, this particular sparse resultant {\bf 
must} vanish due to the fact that our specialization of $F$ results in the 
existence of infinitely many roots of $F$ in $\cT_\wP$. In fact, $F$ 
vanishes on the $1$-dimensional subvariety of $\cT_\wP$ corresponding 
to the left-hand vertical edge of the hexagon $\wP$ (cf.\ section 
\ref{sec:chowtor}). 

Going one dimension higher, consider instead the $3\times 3$ 
system $G$, consisting of the following polynomials:
\begin{eqnarray*}
a_1yz+a_2xz+a_3xy+a_4xyz \\
b_1yz+b_2xz+b_3xy+b_4xyz \\
c_1yz+c_2xz+c_3xy+c_4xyz \\ 
\end{eqnarray*}
Note that the mixed volume bound for this system is again $1$. 

Clearly, $\frac{1}{xyz}G$ is a linear system in $\{\frac{1}{x},
\frac{1}{y},\frac{1}{z}\}$. So by Cramer's rule, we can express 
$x$, $y$, and $z$ as ratios of $3\times 3$ determinants in the 
coefficients. Combining this with the product formula for toric resultants 
\cite{chowprod} (and clearing denominators) we obtain that the sparse 
$u$-resultant of $G$ is precisely\footnote{We
also need the fact that the Pedersen-Sturmfels formula, originally
stated only over $\C$, remains true over a general algebraically
closed field. This is proved in \cite{resvan}.}
\[ [423][143][124]u_0+[123][143][124]u_1+[123][423][124]u_2
+[123][423][143]u_3 \]
where the {\bf bracket} $[ijk]$ \cite{introchow} is the 
$3\times 3$ subdeterminant
\[
\begin{vmatrix} a_i & a_j & a_k \\ b_i & b_j & b_k \\ 
		c_i & c_j & c_k \end{vmatrix}
\]
of the coefficient matrix of $G$. This compactly expressed 
resultant can be thought of as a {\bf semi-mixed} Chow form --- 
a toric resultant of a semi-mixed system \cite{polyhomo}, 
compressed in terms of suitable brackets.  

Now consider the specialization of $G$ to 
\begin{eqnarray*}
yz+xz+2xy+3xyz \\ 
yz+xz+4xy+9xyz \\
yz+xz+8xy+27xyz \\
\end{eqnarray*}
It is then easily verified that $G$ has exactly one 
root in $\cT_P\!\cong\!\Pro^3_K$: 
\[ [x:y:z:w]\!=\![1:-1:0:0] \]
($w$ denoting an extra variable for homogenizing).\footnote{
If $\ch K\!\in\!\{2,3\}$ then $G$ will actually have infinitely 
many roots in $\cT_P$. So let us assume henceforth that 
$\ch K\!\not\in\!\{2,3\}$. (It is easy to construct similar examples when 
$\ch K\!\in\!\{2,3\}$ as well.)} 
More to the point, the sparse $u$-resultant vanishes identically 
for this specialization of $G$, even though $G$ has no zero coefficients.  
Furthermore, one can easily check that the correspondence of (2) (from 
our Main Theorem) does {\bf not} give us the root of $G$ in $\cT_P$. 
(Using $D\!:=\!(\{(0,1,1),(1,1,1)\},
\{(1,0,1),(1,1,1)\},\{(1,1,0),(1,1,1)\})$, the coefficients of 
the $u_i$ in $\cF_A$ suggest the nonsensical root $[0:0:5:21]$.)  

To remedy this, we can use the twisted Chow form $\choo_{A'}(u)$ with 
$$A'\!:=\!\{(0,1,1),(1,0,1),(1,1,0),(1,1,1)\}.$$ In particular, 
when the coefficients of $G$ are unspecialized, 
\[ \choo_{A'}(u)=\begin{vmatrix} a_1 & a_2 & a_3 & a_4\\ 
			    b_1 & b_2 & b_3 & b_4\\
			    c_1 & c_2 & c_3 & c_4\\
			    u_{(0,1,1)} & u_{(1,0,1)} & u_{(1,1,0)} 
			    & u_{(1,1,1)} \\ \end{vmatrix} \]
So under our last specialization, this becomes 
$12u_{(1,0,1)}-12u_{(0,1,1)}$. Thus the coordinates of the sole root of $G$ 
in $\cT_{\conv(A')}\!\cong\!\Pro^3_K$ are precisely $[1:-1:0:0]$. To 
conclude, by virtue of our (compatible) choice of $A'$, the natural embeddings 
of $(\Ks)^3$ into $\cT_P$ and $\cT_{\conv(A')}$ are identical. This is how 
we've recovered our root in $\cT_P$. 

In closing, note that in practice we would never actually 
compute $\cH(u;s)$ --- we would instead recover $\cF_A$ 
via rapid and sophisticated interpolation techniques, e.g., 
\cite{zip}. In particular, our calculations can be sped up tremendously 
with suitably specialized code. 

\subsection{The ``Dense'' Case}
\label{sub:dense}
Our last example illustrates a simple fundamental case.

Suppose $E$ is the $n$-tuple $(d_1\Delta,\ldots,d_n\Delta)$
where $\Delta\!\subset\!\Rn$ is the vertex set of the standard 
$n$-simplex in $\Rn$ and
$d_i\!\in\!\N$ for all $i$. It is then easily verified that the
$n$-tuple $D\!:=\!(\{\bO,d_1\hat{e}_1\},\ldots,\{\bO,d_n\hat{e}_n\})$
is an irreducible fill of $E$ (cf.\ theorem \ref{thm:me}). Letting 
$A\!=\!\Delta$, we see that our polynomial $\cH$ is a variant (over 
a general algebraically closed field) of the original
GCP applied to an $n\!\times\!n$ system of polynomials with degrees
$d_1,\ldots,d_n$ \cite{gcp}. In particular, our $F-sF^*$ has $2n$
$s$-monomials, compared to the $n$ $s$-monomials in Canny's
$(f_1-sx^{d_1}_1,\ldots,f_n-sx^{d_n}_n)$. Note
also that $\conv(A)$ and $P$ are homothetic and $\cT_P\!\cong\!\Pro^n_K$.
Neglecting the extra $s$-monomials, setting $d_i\!=\!1$ for all $i$, and
suitably specializing the coefficients of $g$, we can then recover the usual
characteristic polynomial of a matrix.

\section{Filling}
\label{sec:fill}
Here we briefly recount filling and some
related concepts. Some of the material below is covered at greater length
in \cite{convexapp}. The paper \cite{combiresult}
is also a useful reference but deals more with the sparse resultant than
with filling. The results below form the basis for our combinatorial 
approach to perturbing degenerate polynomial systems.  

Let $\Sn\!\subset\!\Rn$ denote the unit $(n\!-\!1)$-sphere centered at the
origin. For any compact $B\!\subset\!\Rn$ and any $w\!\in\!\Rn$, define
$B^w$ to be the set of $x\!\in\!B$ where the inner-product $x\!\cdot\! w$ is
minimized. (Thus $B^w$ is the intersection of $B$ with its supporting
hyperplane in the direction $w$.) We then define
$E^w\!:=\!(E^w_1,\ldots,E^w_n)$ and
$D\!\cap\!E^w\!:=\!(D_1 \cap E^w_1,\ldots,D_n \cap E^w_n)$.

Recall that the {\bf dimension} of any $B\!\subseteq\!\Rn$, $\dim B$, is
the dimension of the smallest subspace of $\Rn$ containing a translate of
$B$. The following definition is fundamental to our development.
\begin{dfn}
Suppose $C\!:=\!(C_1,\ldots,C_n)$ is an
$n$-tuple of polytopes in $\Rn$ {\bf or} an $n$-tuple of finite subsets of
$\Rn$. We will allow any $C_i$ to be empty and say that a nonempty subset
$J\!\subseteq\![n]$ is {\bf essential} for $C$ (or $C$ {\bf has
essential subset $J$}) $\Longleftrightarrow$ (0) $\supp(C)\!\supseteq\!J$,
(1) $\dim(\sum_{j\in J} C_j)=|J|-1$, and (2) $\dim(\sum_{j\in J'} C_j)\geq
|J'|$ for all nonempty {\bf proper} $J'\!\subsetneqq\! J$.
\end{dfn}

Equivalently, $J$ is essential for $C \Longleftrightarrow$ the
$|J|$-dimensional mixed volume of $(C_j \; | \; j\!\in\!J)$ is $0$ and no
smaller subset of $J$ has this property. Figure \ref{fig:myfirst} below shows 
some simple examples of essential subsets for $C$, for various $C$ in the case
$n\!=\!2$. 

\begin{figure}[bth]
\begin{picture}(410,70)(-35,-5)

\put(-80,-50){
\begin{picture}(133,100)(0,0)
\put(51,80){\circle*{3}} \put(45,85){$C_1$} 
\put(72,80){\circle*{3}} \put(69,85){$C_2$} 
\put(42,50){$\{1\},\{2\}$}
\end{picture}}

\put(24,-50){
\begin{picture}(133,100)(0,0)
\put(51,80){\circle*{3}} \put(45,85){$C_1$}
\put(62,70){\line(1,1){30}}
\put(62,70){\circle*{3}} \put(92,100){\circle*{3}} \put(79,77){$C_2$}
\put(60,50){$\{1\}$}
\end{picture}}

\put(145,-50){
\begin{picture}(133,100)(0,0)
\put(36,70){\line(1,1){30}}
\put(36,70){\circle*{3}} \put(66,100){\circle*{3}} \put(36,85){$C_1$}
\put(62,70){\line(1,1){30}}
\put(62,70){\circle*{3}} \put(92,100){\circle*{3}} \put(79,77){$C_2$}
\put(45,50){$\{1,2\}$}
\end{picture}}

\put(276,-50){
\put(26,100){\line(1,-1){30}}
\put(26,100){\circle*{3}} \put(56,70){\circle*{3}} \put(43,85){$C_1$}
\put(72,70){\line(0,1){30}}
\put(72,70){\circle*{3}} \put(72,100){\circle*{3}} \put(74,85){$C_2$}
\begin{picture}(133,100)(0,0)
\put(45,50){None}
\end{picture}}

\end{picture}

\caption{The essential subsets for 4 different pairs of plane polygons. 
(The segments in the third pair are meant to be parallel.) } 
\label{fig:myfirst}

\end{figure}
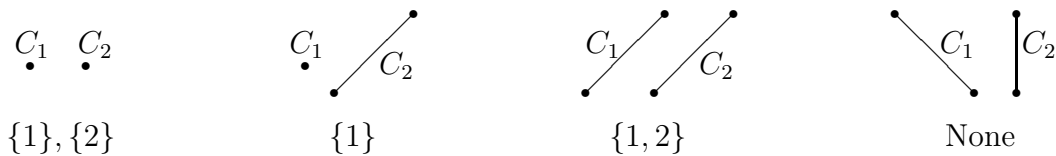

A basic fact about mixed volumes is that $\cM(E)\!=\!0 \Longleftrightarrow
E$ has an essential subset, whenever $\supp(E)\!=\![n]$.
However, there is an even deeper connection between filling and 
essentiality:
\begin{thm} \cite[sec.\ 2.5]{convexapp}
\label{thm:me}
Suppose $D$ and $E$ are $n$-tuples of finite subsets of $\Zn$ such that
$\cM(E)\!>\!0$. Then $D$ fills $E
\Longleftrightarrow$ for all $w\!\in\!\Sn$, $\supp(D\cap E^w)$
contains a subset essential for $E^w$. \qed 
\end{thm}
\begin{rem}
One certainly need not check infinitely many $w$. In fact, we
need only check one $w$ (just pick any inner normal) for each face of
the polytope $\sum^n_{i=1} \conv(E_i)$.
\end{rem}

Filling is closely related to root counting for sparse polynomial 
systems, and this aspect is explored much further in
\cite{convexapp,rojaswang,toricint}. 
We also point out that the computational complexity of finding an irreducible 
fill is an open question. However, for $n\!\leq\!3$, finding 
irreducible fills is quite simple and no harder (asymptotically) than finding 
a convex hull. (Using theorem \ref{thm:me}, this follows as a simple 
geometric exercise.) In any event, the connection between fills and polynomial 
system solving (not to mention specialized resultants) appears
to be new and, we hope, provides added incentive to investigate
filling. Also, even if finding an irreducible fill is too hard, this 
step of our toric GCP construction need 
only be done {\bf once} for a given family of problems, provided
$E$ remains fixed. The situation where the monomial term structure of a
polynomial system remains fixed once and for all, and the coefficients may vary
many thousands of times, actually occurs frequently in many practical 
contexts such as robot control or computational geometry.

\section{Toric Geometry and the Proof of Our Main Theorem} 
\label{sec:chowtor} 
Our notation is a slight variation of that used in 
\cite{tfulton}, and is described at greater length in \cite{toricint}. 
We will assume the reader to be familiar with normal fans of polytopes and  
the construction of a toric variety from a fan or a finite point set 
\cite{tfulton,gkz94}. However, we will at least list our cast of main 
characters:
\begin{dfn} \cite{toricint} 
\label{dfn:toric}
Given any $w\!\in\!\Rn$, we will use the following notation:
\begin{itemize}
\item[$T=$]{The algebraic torus $\Ksn$}
\item[$P^w=$]{The face of $P$ with inner normal $w$} 
\item[$\sigma_w=$]{The closure of the cone generated by the inner normals of 
$P^w$} 
\item[$U_w=$]{The affine chart of $\cT_P$ corresponding to the cone
$\sigma_w$ of $\fan(P)$} 
\item[$L_w=$]{The $\dim(P^w)$-dimensional subspace of $\Rn$ parallel to $P^w$} 
\item[$x_w=$]{The point in $U_w$ corresponding to the semigroup
homomorphism $\sigma^\vee_w\cap\Zn \longrightarrow \{0,1\}$ mapping $p
\mapsto \delta_{w\cdot p,0}$, where $\delta_{ij}$ denotes the Kronecker
delta}  
\item[$O_w=$]{The $T$-orbit of $x_w=$ The $T$-orbit corresponding to $\relint 
P^w$} 
\item[$\cE_P(Q)=$]{The $T$-invariant Weil divisor of $\cT_P$ corresponding to 
a polytope $Q$.} 
\item[$\divisor(f)=$]{The Weil divisor of $\cT_P$ defined by 
a rational function $f$ on $\Ksn$} 
\item[$\cD_P(f,Q)=$]{$\divisor(f)+\cE_P(Q)=$ The toric effective divisor of 
$\cT_P$ corresponding to $(f,Q)$} 
\item[$\cD_P(F,\cP)=$]{The (nonnegative) cycle in the Chow ring of 
$\cT_P$ defined by $\bigcap^k_{i=1} \cD_P(f_i,P_i)$, whenever 
$\cP\!=\!(P_1,\ldots,P_k)$} 
\end{itemize}
\end{dfn}
We will say that a polytope $P$ is {\bf compatible} with $Q$ iff every cone 
of $\fan(Q)$ is a union of cones of $\fan(P)$ \cite{khocompat,tfulton,
toricint}. In particular, whenever $F$ is an $k\times n$ polynomial system with 
support contained in $E$, we will define the {\bf zero set of $\mathbf{F}$ in 
$\mathbf{\cT_P}$} to be the toric cycle $\cD_P(F,\cP)$, 
where $\cP\!:=\!(\conv(E_1),\ldots,\conv(E_k))$. 

The following result will provide some necessary geometric intuition 
for specializing resultants. 
\begin{van} \cite{resvan} 
Suppose $f_i$ is a polynomial over $K$ with support contained in 
$E_i\!\subset\!\Zn$ for all $i\!\in\![n+1]$. Then, provided 
$$\cM(E_1,\ldots,E_{i-1},E_{i+1},\ldots,E_{n+1})\!>\!0$$
for some $i\!\in\![n+1]$, 
\[\res_{\bar{E}}(f_1,\ldots,f_{n+1})\!=\!0 \Longleftrightarrow 
\bigcap\limits^{n+1}_{i=1} \cD_\wP(f_i,\conv(E_i))\!\neq\!
\emptyset,\]
where $\bar{E}\!:=\!(E_1,\ldots,E_{n+1})$ and 
$\wP\!:=\!\sum\limits^{n+1}_{i=1} \conv(E_i)$. \qed 
\end{van}
\begin{rem}
This result provides a geometric analogue, over a {\bf general} 
algebraically closed field, of the product formula for the sparse resultant 
\cite{chowprod}. 
\end{rem}

We will also make frequent use of the natural correspondence 
between the face interiors $\{\relint P^w\}$ and the $T$-orbits 
$\{O_w\}$ \cite{ksz92,tfulton,gkz94}. The following lemma gives a more 
explicit algebraic analogy between the faces of $P$ and the affine 
charts of $\cT_P$.  
\begin{lemma} \cite[Sections 4.2--5.1]{toricint}
\label{lemma:once}
Suppose $F$ is a $k\times n$ polynomial system over $K$ with support 
contained in a $k$-tuple of integral polytopes $\cP\!:=\!(P_1,\ldots,P_k)$ in 
$\Rn$. Assume further that $P$ is a rational polytope in $\Rn$. 
Then the defining ideal in $K[x^e \; | \; e\in\sigma^\vee_w\cap \Zn]$ of 
$U_w\cap \cD_P(F,\cP)$ is $\langle x^{b_i}f_i \; | \;$ for all 
$i\!\in\![k]$ and $b_i\!\in\!\Zn$ such that $b_i+P_i\subset\sigma^\vee_w 
\rangle$. \qed
\end{lemma} 

Lifting (or projecting) from one toric variety to another is an important 
fundamental ideal we will also use. The following lemma follows directly from 
the development of \cite{tfulton}.
\begin{lemma}
\label{lemma:lift}
Suppose $P\!\subset\!\Rn$ is an $n$-dimensional rational polytope, and 
$A$ is either a nonempty finite subset of $\Zn$ {\bf or} a rational 
polytope in $\Rn$. Assume further that $P$ is compatible with $\conv(A)$. Then 
there is a natural (surjective) proper morphism $\fii : \cT_P 
\twoheadrightarrow \cT_A$. In particular, following the notation of this 
section, $\fii(\cD_P(F,\cP))\!=\!\cD_A(F,\cP)$, where the latter cycle is the 
image of $\cD_{\conv(A)}(F,\cP)$ under the natural proper morphism from 
$\cT_{\conv(A)}$ to $\cT_A$. 
\end{lemma} 
\begin{rem}
Recall that $\cT_A$ can be defined as the image of of $\cT_\wP$ under 
the map $\fii$ from our Main Theorem. So, with this understanding, 
there is no ambiguity between our first and second $\fii$. 
\end{rem} 

To conclude our background, we will need the following lemma implying 
that $F^*$ is sufficiently generic in a useful sense. 
\begin{lemma}
\label{lemma:gen}
Suppose $D$ is an irreducible fill of some $n$-tuple $E$. Then for any point 
$v$ lying in any $D_i$, there exists a $w\!\in\!\Rns$ such that $\{i\}$ is the 
unique essential subset of $D^w$ and $D^w_i\!=\!\{v\}$. In particular, 
following the notation of our Main Theorem, $F^*$ has exactly $\cM(\cP)$ roots 
(counting multiplicities) in $\Ksn$. \qed 
\end{lemma}
\noindent
This lemma follows easily from the techniques of \cite{convexapp}, 
particularly section 2.5. 

\subsection{The Proof of Our Main Theorem}
We first note that the well known results on the degree 
of $\res_{\bar{E}}(f_1,\ldots,f_{n+1})$ with respect to the coefficients 
of different $f_i$ \cite{combiresult} remain true over any algebraically 
closed field. This follows easily from the formulation of the resultant for 
a collection of invertible sheafs on a projective variety \cite{gkz94}. 
In particular, the degree of $\cH$ as a polynomial in $s$ should be 
\[ \sum\limits^n_{i=1} \cM(E_1,\ldots,E_{i-1},E_{i+1},\ldots,E_n,A).\] 
Also each coefficient of $\cH(s)$ 
should be a homogeneous polynomial (in the $u_e$) of degree $\cM(E)$. 
These two assertions of course include the opening statement of our Main 
Theorem (on the degree and homogeneity of $\cF_A$), but they will follow 
only upon showing that $\cH$ is not identically $0$. 

To see this, note that lemma \ref{lemma:once} and the Vanishing 
Theorem for Resultants readily imply that the coefficient of 
the {\bf highest} power of $s$ in $\cH$ is precisely 
$\res_{(E,A)}(F^*,g)$. (Simply check the zero set of $F-sF^*$ in 
$\cT_\wP$ at $s\!=\!\infty$.) By lemma \ref{lemma:gen}, and the Vanishing 
Theorem once more, we see that this polynomial in the $u_e$ is not identically 
$0$. So $\cH\!\not\equiv\!0$ and we've finished the simplest part 
of our proof. 

Part (0) of our Main Theorem follows similarly: One need only 
consider the {\bf unspecialized} resultant polynomial $\res_{(E,A)}(F,g)$ 
and observe the terms of degree $0$ in $s$ as we specialize coefficients 
to obtain $F-sF^*$. The statement on the vanishing of $\choo_A$ 
then follows easily from lemma \ref{lemma:lift} (since $\wP$ is 
compatible with $\conv(A)$) and the Vanishing Theorem: We obtain  
that $\fii(Z)$ is positive-dimensional iff $\choo_A$ has infinitely 
many distinct divisors of the form $\sum_{e\in A}c_eu_e$. In particular, 
corollary \ref{dfn:twist} follows as a special case of (0) since $P$ 
compatible with $\conv(A) \Longrightarrow \cT_P\!=\!\cT_\wP$. Note 
that (1), (2), and (3) also follow almost trivially, {\bf provided}
$\choo_A$ is not identically $0$. 

To properly handle the cases of (1), (2), and (3) where we are actually 
working with a non-trivial toric perturbation, let us first construct 
two important toric cycles: Let $\bZ$ be the zero set of $F-sF^*$ 
in $\cT_\wP\times\Pro^1_K$ and $\hZ$ the zero set of $\cH(u;s)$ in 
$\Pa\times\Pro^1_K$. Then it is 
easily observed that $Z\!=\!\bZ\cap (\cT_\wP\times\{0\})$. 
Also, it can be shown \cite[Section 5.1]{toricint} that $\bZ$ contains an 
algebraic curve $C$ (possibly reducible), with surjective projection onto the 
second factor of $\cT_\wP\times\Pro^1_K$, obeying the following property: 
$C\cap(\cT_\wP\times\{0\})$ contains the zero-dimensional part of $Z$, and 
consists of exactly $\cM(E)$ points (counting multiplicities). (In 
fact, $\bZ\cap(\cT_\wP\times\{s_0\})\!=\!C\cap(\cT_\wP\times\{s_0\})$ 
for almost all $s_0\!\in\!\Pro^1_K$.) Furthermore, by 
slightly modifying step (i) of the proof of the toric variety version 
of Bernshtein's theorem \cite[Section 5.1]{toricint}, one can show 
that $C$ intersects every positive-dimensional irreducible component of 
$Z$. (One also needs to use the definition of intersection 
multiplicity of an irreducible component $W$ as the number of curve 
branches intersecting $W$ in a 1-parameter deformation.) 

The proof of the rest of our main theorem will reduce to establishing 
a precise correspondence between the factors of $\cF_A$ and the points 
of $\fii(C\cap(\cT_\wP\times\{0\}))$. To complete this connection, we need 
only observe that $\hZ$ is a very special kind of hypersurface, closely 
related to $\bZ$. 

In particular, if $k$ is the least power of $s$ in $\cH$, observe that
$\hZ$ and the zero set of $\frac{\cH}{s^k}$ in $\Pa\times\Pro^1_K$ differ only 
by the presence of the hyperplane $\Pa\times\{0\}$. The second zero 
set does {\bf not} contain this hyperplane, so let's call the second 
zero set $\dZ$. Then by lemmata \ref{lemma:once} and 
\ref{lemma:lift}, and the Vanishing Theorem for Resultants, $\dim [\bZ\cap 
(\cT_\wP\times\{s_0\})]\!=\!0$ implies the following equivalence: 
$\cH(H_{\fii(\zeta)};s_0)\!=\!0 \Longleftrightarrow 
\zeta\!\in\!\bZ\cap (\cT_\wP\!\times\!\{s_0\})$, where 
$H_p$ is the hyperplane dual to the point $p$.\footnote{So if 
$p\!:=\![p_e \; | \; e\!\in\!A]\!\in\!\Pa$ then $H_p\!:=\!\{[y_e \; | \; 
e\!\in\!A]\!\in\!\Pa \; | \; \sum_{e\in A} p_ey_e\!=\!0\}$.} Now 
$\dim [\bZ\cap (\cT_\wP\times\{\infty\})]\!=\!0$ by construction. Hence 
$\dim [\bZ\cap (\cT_\wP\times\{s_0\})]\!=\!0$ for almost all 
$s_0\!\in\!\Pro^1_K$, by Main Theorem 2 of \cite{toricint}. So 
$\psi^\vee(C)$ is an open subset of $\dZ$, where we define 
$\psi^\vee(C)\!:=\!\{(y,s_0) \; | \; y\!\in\!H_{\fii(\zeta)} \ ; \  
\zeta\!\in\!C\cap(\cT_\wP\times\{s_0\}) \ ; \ s_0\!\in\!\Pro^1_K\}$. Therefore, 
since $\fii$ is a proper map, $\frac{\cH}{s^k}$ must vanish on {\bf all} of 
$\psi^\vee(C)$. In particular, 
\[ \cF_A(u)\!=\!\alpha\cdot\!\!\!\prod\limits_{\zeta\in 
C\cap(\cT_\wP\times\{0\})} \left(\sum\limits_{e\in A} c_{\zeta,e}u_e\right) 
\]
where $\alpha\!\in\!\Ks$, $[c_{\zeta,e} \; | \;
e\!\in\!A]\!:=\!\fii(\zeta)$, and the product counts intersection 
multiplicities.  

Continuing our main proof, (1), (2), and (3) follow immediately 
from our last formula and our preceding observations. \qed  

Note that our algebraic proof avoids the use of limiting 
arguments that were present in \cite{gcp}. Thus our result 
is equally valid when $K$ has positive characteristic. 

\bibliographystyle{amsalpha}

\end{document}